\documentclass[11pt]{article}
\usepackage{paralist} % very flexible & customisable lists (eg. enumerate/itemize, etc.)
\usepackage[utf8]{inputenc}

\usepackage{amssymb,xcolor,xspace}
\usepackage{amsmath,graphicx,comment,bm,%ulem
}
\usepackage[pagebackref]{hyperref}
\usepackage{array}

\usepackage{authblk}

\usepackage{geometry} % to change the page dimensions
\usepackage{fancyhdr} % This should be set AFTER setting up the page geometry
\usepackage{sectsty}
\allsectionsfont{\sffamily\mdseries\upshape} % (See the fntguide.pdf for font help)
\usepackage[nottoc,notlof,notlot]{tocbibind} % Put the bibliography in the ToC
\usepackage[titles]{tocloft} % Alter the style of the Table of Contents

\def\cqfd{\skip10=\parfillskip\parfillskip=0pt
\enspace\hfill\symbolecqfd\par\parfillskip=\skip10\par\medskip}
\def\symbolecqfd{\rlap{$\sqcap$}$\sqcup$}

\newtheorem{theorem}{Theorem}[section]
\newtheorem{proposition}[theorem]{Proposition}
\newtheorem{lemma}[theorem]{Lemma}
\newtheorem{corollary}[theorem]{Corollary}

\newtheorem{pro-fact}[theorem]{Fact}
\newenvironment{fact}{\begin{pro-fact}\rm}{\cqfd\end{pro-fact}}

\newtheorem{pro-example}[theorem]{Example}
\newenvironment{example}{\begin{pro-example}\rm}{\cqfd\end{pro-example}}

\newtheorem{pro-remark}[theorem]{Remark}
\newenvironment{remark}{\begin{pro-remark}\rm}{\cqfd\end{pro-remark}}

\newenvironment{preuve}{\rm \trivlist \item[\hskip \labelsep{\bf
Proof.}]}{\cqfd\endtrivlist}

\def\cqfd{\skip10=\parfillskip\parfillskip=0pt
\enspace\hfill\symbolecqfd\par\parfillskip=\skip10\par\medskip}
\def\symbolecqfd{\rlap{$\sqcap$}$\sqcup$}
\def\proof{\begin{preuve}}
\def\eop{\end{preuve}}

\let\phi\varphi

\def\inv{^{-1}}
\let\epsilon\varepsilon
\def \calA {\mathcal{A}}

\def \O {\mathcal{O}}

\def\N{\mathbb{N}}

\def\pr{\textsf{pr}}

\def\probap{\mathfrak{p}}
\def\probaq{\mathfrak{q}}

\def\Maa{M_{a,a}}
\def\Mab{M_{a,b}}

\def\mf{\textsf{mf}}

\def\Aut{\textsf{Aut}}

\def\path{\textsf{path}}

\def \calG {\mathcal{G}}

\def\algoPPP{\mathcal{PPP}}
\def\algoMP{\mathcal{MP}}
\def\algoMPd{\mathcal{MP}_d}
\def\algoPrim{\mathcal{P}}
\def\algoShp{\mathcal{S}}
\def\algoCR{\mathcal{CR}}
\def\algoRPd{\mathcal{RP}_d}

\title{The central tree property and algorithmic problems on subgroups of free groups}

\author[1]{Mallika Roy\thanks{mallikaroy75@gmail.com}}
\author[2]{Enric Ventura\thanks{enric.ventura@upc.edu}}
\author[3,4]{Pascal Weil\thanks{pascal.weil@cnrs.fr}}

\affil[1]{Departamento de Matem\'aticas, UPV/EHU, Bilbao, Spain}
\affil[2]{Department de Matem\`atiques, and Institut de Matem\`{a}tiques de la UPC-BarcelonaTech, Univ. Polit\`ecnica de Catalunya, Catalonia}
\affil[3]{CNRS, ReLaX, IRL 2000, Siruseri, India}
\affil[4]{Univ. Bordeaux, CNRS, Bordeaux INP, LaBRI, UMR 5800, F-33400 Talence, France}
\date{} % delete this line to display the current date

\begin{document}

\date{}

\maketitle

\begin{abstract}
We study the average case complexity of the uniform membership problem for subgroups of free groups, and we show that it is orders of magnitude smaller than the worst case complexity of the best known algorithms. This applies to subgroups given by a fixed number of generators as well as to subgroups given by an exponential number of generators. The main idea behind this result is to exploit a generic property of tuples of words, called the central tree property.

An application is given to the average case complexity of the relative primitivity problem, using Shpilrain's recent algorithm to decide primitivity, whose average case complexity is a constant depending only on the rank of the ambient free group.
\end{abstract}

%%%%%%%%%%%%%%%%%%
\section{Introduction}

Algorithmic problems have been prominent in the theory of infinite groups at least since Dehn formulated the word problem for finitely presented groups \cite{1912:Dehn}: a finite group presentation $\langle A \mid R \rangle$ being fixed, the word problem asks whether a given element of the free group on $A$ --- seen as a reduced word on the alphabet $A \cup A\inv$ --- is equal to the identity in the group presented by $\langle A \mid R \rangle$.

A problem is called \emph{decidable} if one can exhibit an algorithm that solves it. It is well known that the word problem is undecidable for certain finite group presentations (Novikov \cite{1955:Novikov}). For decidable problems, it is natural to try to evaluate the complexity of an algorithm solving them, namely the amount of resources (time or space) required to run the algorithm, as a function of the size of the input.

The most common complexity evaluation for an algorithm $\calA$, is the \emph{worst case complexity}, which measures the maximum time required to run $\calA$ on an input of size $n$. In certain cases, it may be relevant to consider the \emph{generic complexity} of $\calA$: $\calA$ has generic complexity at most $f(n)$ if the ratio of inputs of size $n$ on which $\calA$ requires time at most $f(n)$ tends to 1 as $n$ tends to infinity. This notion of complexity recognizes that the instances that are hard for $\calA$ (the witnesses of the worst case complexity) may be few but it does not attempt to quantify the time required on the instances in a vanishing set.

Here we will be concerned with the more precise \emph{average case complexity}, namely the expected time required to run $\calA$ on size $n$ instances taken uniformly at random. This measure of complexity takes into account the resources needed for every input.

The specific problems we consider in this paper are the uniform membership problem and the relative primitivity problem in finite rank free groups. Recall that an element $w$ of $F(A)$, the free group on $A$, is \emph{primitive} if $F(A)$ admits a free basis containing $w$. The \emph{Uniform Membership Problem} (resp., the \emph{Relative Primitivity Problem}) is the following: given elements $w_0, w_1, \dots, w_k$ of $F(A)$, decide whether $w_0$ belongs to (resp., is primitive in) the subgroup $H$ of $F(A)$ generated by $w_1,\dots, w_k$. In this paper, the length $k$ of the tuple $(w_1,\ldots, w_k)$ is not fixed, and the parameters we consider to gauge the size of an instance are $n=\max\{|w_i| \mid 1\le i\le k\}$, $k$ as a function of $n$, and $m=|w_0|$. In particular, the total length of an input $(w_0,w_1,\ldots,w_k)$ is at most $kn+m$.

It is well known that both these problems are decidable, and can be solved in polynomial time worst case complexity (see \cite{1983:Stallings}, \cite[Fact 3.6]{2007:RoigVenturaWeil}, and Sections~\ref{sec: gmp}, \ref{sec: primitivity} below).
The most efficient solution (again, from the point of view of the worst case complexity) of the Uniform Membership Problem uses the concept of Stallings graph of a subgroup, a finite $A$-labeled graph uniquely associated with a finitely generated subgroup of $F(A)$, which can be easily computed and then gives a linear time solution for the uniform membership problem. More precisely, the Stallings graph of $H$ has at most $kn$ vertices and it is computed in time $\O(kn\log^*(kn))$ \cite{2006:Touikan}. After this computation, deciding whether $w_0\in H$ is done in linear time in $m$.

Our main result is an algorithm solving the Uniform Membership Problem, whose average case complexity is asymptotically a little $o$ of the worst case complexity described above, at least when the number $k$ grows at most polynomially with $n$. It is notable that the dependance of our algorithm's expected performance on $m$ (the length of the word $w_0$ to be tested) is extremely low.

A specific instance of our main result shows, for instance, that if $k$ is a constant, then the Uniform Membership Problem can be solved in expected time $\O(\log n + mn^{-\log(2r-1)})$, where $r=|A|$ is taken to be constant.

The fundamental ingredient in this result is the so-called \emph{central tree property} (ctp) for a tuple $(w_1,\dots,w_k)$. This property, formally introduced in \cite{2016:BassinoNicaudWeilCM,2020:BassinoNicaudWeil}, holds when the $w_i$ and their inverses have little initial cancellation (that is: they have short common prefixes). This property turns out to hold with high probability and, when it holds, solving the Uniform Membership Problem is considerably simpler than in the general case.

We then apply our result to the Relative Primitivity Problem: we give an algorithm solving it whose average case complexity is much lower than its worst case complexity. Here, we use an algorithm recently proposed by Shpilrain \cite{2023:Shpilrain} to solve the primitivity problem (deciding whether a given word $w_0$ is primitive in $F(A)$). Shpilrain's algorithm has the remarkable property of having constant average case complexity; that is, its expected time does not depend on the length of $w_0$. As it turns out, this constant average case complexity depends on the rank of the ambient free group and this is important in the context of the Relative Primitivity Problem, where we need to test for primitivity in the subgroup generated by $w_1,\dots,w_k$, whose rank may be as large as $k$.

The paper is organized as follows. Section~\ref{sec: prelim} briefly discusses the fundamental notions on subgroups of free groups which we will use, especially the notions of the Stallings graph and the growth function of a subgroup, as well as the notion of average case complexity and the computational model which we rely upon.

Section~\ref{sec: ctp} is dedicated to the central tree property, applied to a tuple $\vec w=(w_1,\dots,w_k)$: its definition, its consequences in terms of the rank and the growth function of the subgroup generated by $\vec w$, and the probability that it holds (in terms of the parameters $k$ and $n$).

Our main result on the average case complexity of the Uniform Membership Problem is presented in Section~\ref{sec: gmp}. Finally we discuss Shpilrain's constant average case complexity algorithm for the primitivity problem, and our application to the Relative Primitivity Problem in Section~\ref{sec: primitivity}.

%%%%%%%%%%%%%%%%%%
\section{Preliminaries}\label{sec: prelim}

%%%%%%%%%%%%%%%%%%
\subsection{Subgroups of free groups}\label{sec: subgroups of free groups}

Throughout the paper, $A$ is a finite non-empty set, called an \emph{alphabet}. We say that a directed graph $\Gamma$ is an {$A$-graph} if its edges are labeled with elements of $A$.  A \emph{rooted} $A$-graph is a pair $(\Gamma,v)$ where $\Gamma$ is a finite $A$-graph and $v$ is a vertex of $\Gamma$. Finally, we say that a rooted $A$-graph $(\Gamma,v)$ is \emph{reduced} if $\Gamma$ is finite and connected, distinct edges with the same start (resp., end) vertex always have distinct labels, and every vertex except possibly the root $v$, is incident to at least two edges.

The set $\tilde A = \{a, a\inv \mid a\in A\}$ (with cardinality $2|A|$), is called the (\emph{symmetrized}) \emph{alphabet}, its elements are called \emph{letters}. We denote by $\tilde A^*$ the set of \emph{words} on $\tilde A$, that is, of finite sequences of letters. We also denote by $A^*$ the set of words using only letters in $A$. A word in $\tilde A^*$ is said to be \emph{reduced} if no letter $a\in A$ is immediately preceded or followed by the letter $a\inv$.

For convenience, if $A = \{a_1,\dots, a_k\}$, we let $a_{-i} = a_i\inv$ for $1\le i\le k$, so that $\tilde A = \{a_i \mid -k \le i \le k,\ i\ne 0\}$.

Suppose that $w=x_1\cdots x_m$ is a word in $F(A)$ (with each $x_i \in \tilde A$) and that $p,q$ are vertices of a reduced $A$-graph $\Gamma$. We say that $w$ \emph{labels a path in $\Gamma$ from $p$ to $q$} if there exists a sequence of vertices $p_0 =p, p_1, \dots, p_m =q$ such that, for every $1\le i\le m$, $\Gamma$ has an $x_i$-labeled edge from $p_{i-1}$ to $p_i$ if $x_i \in A$, and an $x_i\inv$-labeled edge from $p_i$ to $p_{i-1}$ if $x_i\inv \in A$. If the start and end vertices of the path are equal (that is, if $p = q$), we say that the path is a \emph{circuit at} $p$.

The free group on $A$ is written $F(A)$, and we identify it with the set of reduced words in $\tilde A^*$. It is well known that every subgroup of $F(A)$ is free \cite{1918:Nielsen}. It is also well-known that every finitely generated subgroup $H$ of $F(A)$ can be associated with a uniquely defined reduced rooted $A$-graph $(\Gamma(H),1)$, called the \emph{Stallings graph} of $H$, with the following property: a reduced word is in $H$ if and only if it labels a circuit in $\Gamma(H)$ at vertex $1$. We refer the reader to the seminal works of Serre \cite{1980:Serre} and Stallings \cite{1983:Stallings} who introduced this combinatorial tool, and to \cite{2022:DelgadoVentura,2002:KapovichMyasnikov,2001:MargolisSapirWeil,2007:MiasnikovVenturaWeil,2007:RoigVenturaWeil} for some of its many applications.

Of particular interest for this paper are the following facts.
\begin{itemize}
\item Given $\vec w=(w_1,\dots,w_k)$ a tuple of reduced words in $F(A)$, one can effectively compute the Stallings graph $(\Gamma(H),1)$ of the subgroup $H$ generated by the $w_i$ and this graph has at most $n = \sum_i|w_i|$ vertices. Touikan \cite{2006:Touikan} showed that $(\Gamma(H),1)$ can be computed in time $\O(n\log^*n)$. Recall that $\log^*(n)$ is the least integer $k$ such that the $k$-th iterate of the logarithmic function yields a result less than or equal to 1: that is, such that $\log^{(k)}(n)\leq 1 <\log^{(k-1)}(n)$. Equivalently, let $b_1=2$ and $b_{k+1}=2^{b_k}$. Then $\log^*n=1$ if $n\leq b_1$ and $\log^*n=k$ if $b_{k-1}<n\leq b_k$.

\item Once $\Gamma(H)$ is constructed, deciding whether a word $w_0\in F(A)$ is an element of $H$ is done by checking whether $w_0$ labels a circuit in $\Gamma(H)$ at vertex $1$. This can be done in time $\O(|w_0|)$.

\item Let $V$ and $E$ be the sets of vertices and edges of $\Gamma(H)$, respectively. Let $T$ be a spanning tree of $\Gamma(H)$ (that is, a subgraph of $\Gamma(H)$ which is a tree and contains every vertex of $\Gamma(H)$). Let $E_T$ be the set of edges of $T$: then we have $|E_T|=|V|-1$. For each vertex $p$ of $\Gamma(H)$, let $u(p)$ be the only reduced word which labels a path in $T$ from the root vertex 1 to vertex $p$ (with $u(1)$ the empty word). For each edge $e$ of $\Gamma(H)$ which is not in $E_T$, say, $e$ is an edge from vertex $p_e$ to vertex $q_e$ with label $a(e)\in A$, let $b(e)=u(p_e)a(e)u(q_e)\inv$. Then $b(e)$ is a reduced word in $F(A)$ which labels a circuit at 1, so that $b(e) \in H$. Moreover, the set $\{b(e) \mid e\in E\setminus E_T\}$ is a basis of $H$.

\item If $T$ is a fixed spanning tree of $\Gamma(H)$ and $B$ is the corresponding basis of $H$, in bijection with $E \setminus E_T$, the expression of an element of $H$ in that basis is obtained as follows. Given a reduced word $w$ in $H$, consider the circuit at 1 labeled by $H$ and the sequence of edges not in $T$ traveled by this circuit, say $e_1^{\epsilon_1},\ldots ,e_h^{\epsilon_h}$, where the $e_i$ are in $E \setminus E_T$, $\epsilon_i =1$ if $e_i$ is traversed in the direct sense, and $\epsilon_i =-1$ if $e_i$ is traversed backwards. Then $w=b(e_1)^{\epsilon_1}\cdots b(e_h)^{\epsilon_h}$.

\item One can compute a spanning tree of $\Gamma(H)$ in time $\O(|V|+|E|)$ (by classical depth-first search). Since $|E| \le |V|\,|A|$, it follows that one can compute a spanning tree of $\Gamma(H)$, and therefore a basis of $H$, in time $\O(|V|\,|A|)$.

\item The subgroup $H$ has finite index if and only if every vertex of $\Gamma(H)$ is the origin (and so, the terminus) of an $a$-labeled edge for every letter $a\in A$ or, equivalently, if $|E|=|V|\,|A|$. In that case, the index of $H$ is $|V|$.
\end{itemize}

%%%%%%%%%%%%%%%%%%
\subsection{Growth modulus}\label{sec: growth modulus}

The \emph{growth function} of a set $L$ of words over an alphabet $A$ is the function $s_L(n)$ counting the words of length $n$ in $L$. The following result belongs to the folklore of combinatorial automata theory.

\begin{fact}\label{fact: growth function of regular languages}
If $L$ is a regular language, then its growth function (restricted to its support, which is the complement of an ultimately periodic sequence) is asymptotically equivalent to an expression of the form $C n^k \lambda^n$, with $C > 0$, $k \in \N$ and $\lambda\ge 1$.
\end{fact}

The real number $\lambda$ in Fact~\ref{fact: growth function of regular languages} is called the \emph{growth modulus} of $L$.

\proof A result usually attributed to Chomsky and Sch\"{u}tzenberger \cite{1963:ChomskySchutzenberger} (see \cite[Proposition I.3]{2009:FlajoletSedgewick} for a quick proof) states that the generating function of $L$, $S_L(z) = \sum_{n\ge 0} s_L(n)z^n$, is a rational fraction. The partial fraction decomposition (over the reals) of this rational fraction yields the announced asymptotic equivalent for the coefficients $s_L(n)$ of $S_L(z)$.
\eop

We record the following elementary remark.

\begin{remark}\label{rk: length at most n}
If a language $L$ has at most $C\lambda^n$ words of length $n$ ($C>0$, $\lambda>1$), then the number of its words of length less than or equal to $n$ is at most $C\frac\lambda{\lambda-1}\lambda^n$, which is $\Theta(\lambda^n)$.
\end{remark}

The Perron-Frobenius theorem (see, e.g., \cite[Theorem 8.4.4] {1990:HornJohnson}) yields a more precise characterization of the growth modulus of a regular language. If $\calA$ is a finite state automaton over alphabet $A$ with state set $Q$, we denote by $G_\calA$ the graph with vertex set $Q$ and with an edge from vertex $p$ to vertex $q$ for every letter $a\in A$ labeling a transition from $p$ to $q$. We also let $M_\calA$ be the associated incidence $(Q\times Q)$-matrix. Note that $M_\calA$ has only non-negative integer coefficients. We say that $\calA$ is \emph{irreducible} if $M_\calA$ is, that is, if $G_\calA$ is strongly connected.

The \emph{period} of $\calA$ is defined to be the largest positive integer $d$ such that $Q$ can be partitioned as $Q=Q_1\sqcup\dots\sqcup Q_d$ in such a way that every transition from a state in $Q_i$ leads to a state in $Q_{i+1}$ (indices are taken modulo $d$). For instance, the minimal automaton of the language of words of length a multiple of $d$, has period $d$. Finally, we say that $\calA$ is aperiodic if its period is 1.

The following result also belongs to the folklore, see \cite[Proposition V.7]{2009:FlajoletSedgewick}.

\begin{proposition}\label{prop: perron-frobenius}
Let $L$ be a regular language, accepted by a deterministic finite state automaton $\calA$ which is irreducible and aperiodic. The growth modulus of $L$ is equal to the dominant eigenvalue of the transition matrix of $\calA$.
\end{proposition}

The growth modulus of the free group $F(A)$ is easily computed.

\begin{example}\label{ex: growth modulus of F(A)}
It is clear that, for every integer $n\geq 1$, the number $R_n$ of reduced words of length $n$ is $2r(2r-1)^{n-1}$, where $r = |A|$.

Recall that we identify $F(A)$ with the set of reduced words over the alphabet $\tilde A$. Thus the growth modulus of $F(A)$ is $2r-1$.
\end{example}

%%%%%%%%%%%%%%%%%%
\subsection{Algorithmic problems: average case complexity}\label{sec: algorithmic problems}

In evaluating the complexity of an algorithm, one needs to specify the model of computation, and the input space. The input space is usually equipped with a notion of size (a positive integer), such that there are finitely many inputs of any given size. It is important, in particular, to make it clear which parameters of the problem are taken to be constants. Unless otherwise indicated, we consider in this paper that the rank $r = |A|$ of the ambient free group is a constant.

The model of computation we adopt is the standard RAM model. Concretely, this means that if a word $w$ is part of the input, it takes unit time to move the reading head to a position $i$ (where $i$ has been computed before), and unit time to read the letter of $w$ in position $i$. Arithmetic operations on integers (addition, multiplication) are also considered as taking unit time.

\begin{remark}\label{rk: length of words in input}
In the standard RAM model, figuring out the length of an input word $w$ (in order, for instance, to read its last letter) takes time $\O(\log|w|)$, using an instance of the so-called exponentiation search method \cite{1976:BentleyYao}: one reads letters in positions 1, 2, 4, 8, 16, etc, until one exceeds the length of the word after, say, $b$ steps At that point, we know that $2^{b-1} \le |w| < 2^b$, that is, we know the leading bit of the binary expansion of $|w|$. The next bits are established by a classical dichotomy method. Concretely, probing position $c=2^b+2^{b-1}$ allows us to know the second bit: $1$ if position $c$ is still in the word, $0$ if $c$ exceeds its length. This is repeated for the successive bits of the binary expansion of $|w|$.

Since we are going to work with algorithms with very low, even constant, average case complexity, we do not want to have to add the logarithmic time needed to compute the length of input words, and will therefore accompany every input word with its length.
\end{remark}

\begin{remark}\label{rk: bitcost model}
When handling ``very large'' integers, or words on a ``very large'' alphabet, it may be more appropriate to use the bitcost model: adding integers takes time linear in the length of their binary representations (that is: in their logarithm) and reading or comparing letters from a large alphabet $A$ takes time proportional to $\log |A|$ (since each letter can be encoded in a bit string of length $\lceil\log|A|\rceil$).
\end{remark}

The \emph{worst case complexity of an algorithm} $\calA$ is the function $f(n)$, defined on $\N$, which accounts for the maximum time required to run Algorithm $\calA$ on an input of size $n$. If a distribution is specified on the set of size $n$ inputs (in this paper: the uniform distribution), the \emph{average case complexity} of $\calA$ is the function $g(n)$ which computes the expected time required to run $\calA$ on inputs of size $n$. Such complexity functions are usually considered up to asymptotic equivalence.

The average case complexity of an algorithm is obviously bounded above by the worst case complexity, and it may sometimes be much lower. The usual idea in discussing average case complexity is to distinguish, within the input space, between a subset of high probability where the algorithm performs very fast, and its low-probability complement, containing all the hard instances (those which witness the worst case complexity).

Finally, the worst case (resp., average case) \emph{complexity of a problem} is the lowest worst case (resp., average case) complexity of an algorithm solving this problem.

A well-known example which will be useful in the sequel, is the \emph{Proper Prefix Problem} (PPP) on alphabet $A$: given two words $u,v$ on alphabet $A$, decide whether $u$ is a proper prefix of $v$ (that is, $v = uu'$ for some non-empty word $u'$). We also consider the \emph{Prefix Problem} (PP) --- deciding whether $u$ is a prefix of $v$ --- and the \emph{Equality Problem} (EqP)  --- deciding whether $u=v$. The following observation is elementary: it is just an instance of the fact that the expected value of a geometric distribution is constant.

\begin{lemma}\label{lm: constant time prefix}
Let $A$ be a finite alphabet with $|A|\geq 2$ and suppose that each set $A^n$ (the words of length $n$) is equipped with the uniform distribution. Problems PPP, PP and EqP can be solved in constant expected time.

The same holds for the set $R_n$ of (reduced) words of length $n$ in $F(A)$.
\end{lemma}

\proof We prove the statement for PPP. The proofs for PP and EqP are entirely similar. Here is a simple and natural algorithm solving the PPP.

\paragraph{Algorithm $\algoPPP$}
On input $u,v$, read $u$ and $v$ from left to right one letter at a time, comparing each letter of $u$ with the corresponding letter of $v$, and stopping because (i) it detected a difference (that is: the $i$-th letters of $u$ and $v$ are different for some $i$); (ii) it reached the end of $u$ but not the end of $v$; or (iii) it reached the end of $v$.

\medskip

It is clear that $u$ is a proper prefix of $v$ in Case (ii), and not a proper prefix of $v$ in every other case. That is: $\algoPPP$ solves the PPP. We now show that it has constant average case complexity. In effect, this is due to the fact that, with high probability, we will detect a difference between the words $u$ and $v$ without having to read either word to its end.

Let $u_i$ (resp., $v_i$) denote the $i$-th letter in $u$ (resp., $v$). At each step, verifying whether $u_i = v_i$ is done in constant time, and it is the case with probability $p = \frac1{|A|}$. Therefore the probability that the algorithm stops after exactly $k$ steps (with $k\leq |u|,|v|$) is $p^{k-1}(1-p)$. Detecting whether we reached the end of $u$ or $v$ is also done in constant time. It follows that the average case complexity of $\algoPPP$ is bounded above, up to a multiplicative constant, by
 $$
1+p+\cdots +p^{|u|-1} \leq \frac1{1-p}=\frac{|A|}{|A|-1} \leq 2,
 $$
and this concludes the proof relative to $A^*$.

To transfer this result to $F(A)$, one can for instance rewrite the words in $F(A)$ as follows. Let $X = \{x_1,\dots,x_{2r-1}\}$ and $\tilde A = \{a_{-r},\dots,a_{-1},a_1,\dots a_r\}$ be ordered in the natural way. For each $b\in \tilde A$, we let $\zeta_b$ be the order isomorphism from $\tilde A \setminus\{b\}$ to $X$.

Given a word $w=b_1\cdots b_n \in F(A)$ of length $n\ge 2$, we let $\tau(w)=b_1c_2\cdots c_n$ be the word in $\tilde AX^{n-1}$ where $c_{i+1}=\zeta_{b_i\inv}(b_{i+1})$ for every $1\leq i<n$. It is clear that $\tau$ is a bijection between $R_n$ and $\tilde AX^{n-1}$ and hence it preserves the uniform distribution.

The result follows if we modify Algorithm $\algoPPP$ as follows: on input $u$ and $v$, the algorithm first compares $u_1$ and $v_1$, and then compares $\zeta_{u_i\inv}(u_{i+1})$ and $\zeta_{v_i\inv}(v_{i+1})$ until a difference is detected.
\eop

\begin{remark}\label{rk: prefix over large alphabets}
Lemma~\ref{lm: constant time prefix} establishes that the expected running time of $\algoPPP$ is bounded by a constant. This relies on our choice of the RAM model of computation, according to which it takes constant time to compare two letters in $A$. As mentioned in Remark~\ref{rk: bitcost model}, if $A$ is very large, it may take non-trivial time to perform such a comparison (namely time $\O(\log|A|)$) and, in that case, the expected running time of $\algoPPP$ is $\O(\log|A|)$.
\end{remark}

%%%%%%%%%%%%%%%%%%
\section{The central tree property: a generic property of tuples of words}\label{sec: ctp}

If $d$ is a positive integer, we say that the $k$-tuple $\vec w=(w_1,\dots,w_k)$ of words in $F(A)$ has the \emph{central tree property} \emph{of depth $d$} (the \emph{$d$-ctp} for short) if the $w_i$ have length greater than $2d$ and the prefixes of length $d$ of the $w_i$ and the $w_i\inv$ are pairwise distinct. We also say that \emph{$\vec w$ has the ctp} if it has the $d$-ctp for some $d<\frac12\min\{|w_i| \mid 1\leq i\leq k\}$. The central tree property was formally introduced in~\cite{2016:BassinoNicaudWeilCM} (see also~\cite{2020:BassinoNicaudWeil}) but it was implicit in the literature, especially on the (exponential genericity of the) small cancellation property, since the ctp can be viewed as a \emph{small initial cancellation property}.

Let $\vec w=(w_1,\ldots,w_k)$ be a $k$-tuple of words in $F(A)$. For convenience, we let $\min|\vec w| = \min\{|w_i| \mid 1\le i\le k\}$ and $\max|\vec w| = \max\{|w_i| \mid 1\le i\le k\}$, and we write $w_{-i}$ for $w_i\inv$ ($1\le i\le k$).

Suppose that $\vec w$ has the $d$-ctp. Then we let $\pr_i$ be the length $d$ prefix of $w_{-i}$ and $\mf_d(w_i)$ be the \emph{middle factor} of $w_i$ of length $|w_i|-2d$. In particular, we have $\mf_d(w_{-i}) = \mf_d(w_i)\inv$ and
\begin{equation}\label{eq: factoring w_i}
w_i =\pr_{-i} \cdot \mf_d(w_i) \cdot \pr_i\inv,
\end{equation}
for $-k\leq i\leq k$, $i\neq 0$.

Denote by $L(\vec w)$ the set of all the $\pr_{i}$. By definition of the $d$-ctp, $|L(\vec w)| = 2k$. Let also $\Gamma_d(\vec w)$ be the tree of prefixes of the $w_i$ and $w_i\inv$ (rooted at the empty word): this is the graph with vertices all the prefixes of the words $w_i$ and $w_i\inv$ ($1\leq i\leq k$), including the empty word, and with an edge from a word $v$ to a word $w$ exactly if $w=va$ ($a\in \tilde A$). We identify the words in $L(\vec w)$ with the corresponding leaves of $\Gamma_d(\vec w)$. If $H=\langle\vec w\rangle$ (i.e., $H=\langle w_1,\dots,w_k \rangle$), then $\Gamma(H)$ consists of the \emph{central tree} $\Gamma_d(\vec w)$, together with $k$ disjoint paths: for every $1\leq i\leq k$, there is such a path from vertex $\pr_{-i}$ to vertex $\pr_i$, labeled $\mf_d(w_i)$. In view of Equation~\eqref{eq: factoring w_i}, the word $w_i$ labels a circuit at the root in $\Gamma(H)$, going first through $\Gamma_d(\vec w)$ to the leaf $\pr_{-i}$, and returning to the root through the leaf $\pr_i$.

\begin{remark}\label{rk: O(kd)}
It follows directly from the definition that one can decide whether a $k$-tuple $\vec w$ has the $d$-ctp, and construct $\Gamma_d(\vec w)$, in time $\O(kd)$.
\end{remark}

We record the following property of $k$-tuples with the ctp, the first of which is~\cite[Lemma 1.2]{2016:BassinoNicaudWeilCM}.

\begin{proposition}\label{prop: ctp basis}
Let $d\geq 1$, let $\vec w$ be a tuple of words in $F(A)$ with the $d$-ctp and let $H=\langle \vec w\rangle$. Then $H$ has infinite index and $\vec w$ is a basis of $H$.

Let $\mu =\min|\vec w|$ and $\nu =\max|\vec w|$. If $w_0$ is a word in $F(A)$ which belongs to $H$, then the length $\ell$ of the expression of $w_0$ in the basis $\vec w$ satisfies $(\mu-2d)\ell \leq |w_0|\leq \nu\ell$.
\end{proposition}

\proof The infinite index property is immediately verified, since $\Gamma(H)$ has vertices of degree 2, namely the leaves of $\Gamma_d(\vec w)$, see Section~\ref{sec: subgroups of free groups}.

Let $X=\{x_1,\dots,x_k\}$ be a $k$-letter alphabet and let $\phi\colon F(X)\to F(A)$ be the morphism given by $\phi(x_i)=w_i$. It is obvious that the image of $\phi$ is $H$. Recall that we let $x_{-i}=x_i\inv$ for every $1\leq i\leq k$. Let $x=x_{i_1}\cdots x_{i_\ell}$ be a non-empty word in $F(X)$. Then $\phi(x)$ is the word obtained by reducing $w_{i_1}\cdots w_{i_\ell}$. Because of the ctp, reduction occurs only in segments of length at most $2d$ around the boundary between $w_{i_h}$ and $w_{i_{h+1}}$, for each $1\leq h<\ell$. In particular, $\ell\mu -2d(\ell-1) \leq |\phi(x)|\leq \ell\nu$.

It follows that $\phi(x) \neq 1$, and hence $\phi$ is injective and $\vec w$ is a basis of $H$.
\eop

It is well-known that an infinite index subgroup $H$ of $F(A)$ has growth modulus smaller than $2r-1$. In the case $H$ is generated by a tuple with the ctp, its growth modulus is greatly constrained.

\begin{proposition}\label{prop: growth under ctp}
Let $\vec w$ be a $k$-tuple of words with the $d$-ctp and let $\mu =\min|\vec w|$. Let $H=\langle\vec w\rangle$. Then the growth modulus of $H$ is at most $(2k-1)^{\frac{1}{\mu-2d}}$.
\end{proposition}

\proof
Let $X$ and $\phi$ be as in the proof of Proposition~\ref{prop: ctp basis}, let $x\in F(X)$ and $w=\phi(x)\in H$. Then $(\mu-2d)|x|\leq |w|$. If $|w|=m$, then $|x|\leq \frac{m}{\mu-2d}$. In particular, the set of words of $H$ of length $m$ is contained in the $\phi$-image of the words in $F(X)$ of length at most $\frac{m}{\mu-2d}$. This set has cardinality $\Theta((2k-1)^{\frac m{\mu-2d}})$, see Remark~\ref{rk: length at most n}. The announced inequality follows since $\phi$ is a bijection between $F(X)$ and $H$.
\eop

We also record the following fact, which is elementarily verified and is well known (see, \textit{e.g.}, \cite{1996:ArzhantsevaOlshanskii}).

\begin{proposition}\label{prop: long words in a tuple}
Let $k \ge 1$. The probability, for a $k$-tuple $\vec w$ of words in $F(A)$ of length at most $n$, to satisfy $\min|\vec w| \le n/2$ is $\O\left(k(2r-1)^{-\frac{n}2}\right)$.
\end{proposition}

\proof The set of words of length $h$ is $2r(2r-1)^{h-1}$, so the set of words of length at most $h$ is
 $$
1+2r+2r(2r-1)+\cdots +2r(2r-1)^{h-1}=\frac{r}{r-1} (2r-1)^h -\frac{1}{r-1}.
 $$
The probability that a word in $F(A)$ of length at most $n$, actually has minimal length at most $n/2$ is therefore asymptotically equivalent to $C\,(2r-1)^{n/2-n}=C\, (2r-1)^{-n/2}$ for some constant $C$, and the result follows.
\eop

We are interested in the $d(n)$-ctp, where $d(n)$ is an increasing function of $n$. The following statement is derived from \cite{2016:BassinoNicaudWeilCM}.

\begin{proposition}\label{prop: failing ctp}
Let $r=|A|\geq 2$ and let $k\ge 2$ be an integer. Let $d(n)$ be a non-decreasing function of $n$ such that $d(n)<n/2$. A random $k$-tuple of words in $F(A)$ of length at most $n$ fails the $d(n)$-ctp with probability $\O(k^2(2r-1)^{-d(n/2)})$.
\end{proposition}

\proof Let $\eta_n =k^2(2r-1)^{-d(n/2)}$. It is shown in~\cite[proof of Prop.~3.17]{2016:BassinoNicaudWeilCM}\footnote{Proposition 3.17 in \cite{2016:BassinoNicaudWeilCM} is formulated in the case of a so-called prefix-heavy distribution on words of fixed length. This is the case for the distribution used here, namely the uniform distribution on reduced words of fixed length. The parameters $C$ and $\alpha$ in~\cite[Proposition 3.17]{2016:BassinoNicaudWeilCM} are, respectively 1 and $(2r-1)\inv$, see~\cite[Example 3.2]{2016:BassinoNicaudWeilCM}.} that the probability that a $k$-tuple $\vec w=(w_1,\ldots,w_k)$ of words in $F(A)$ of length at most $n$ fails to have the $d(n)$-ctp is bounded above by the sum of $5\eta_n$ and the probability that $k^2 (2r-1)^{-d(\min|\vec w|)}>\eta_n$.

We note that, if $\min|\vec w|>n/2$, then for all $n$, we have $k^2(2r-1)^{-d(\min|\vec w|)} \leq k^2(2r-1)^{-d(n/2)}=\eta_n$.

Therefore, the probability that $\vec w$ fails to have the $d(n)$-ctp is bounded above by the sum of the probability that $\min|\vec w|\leq n/2$, which is $\O\left( k(2r-1)^{-\frac{n}2} \right)$ by Proposition~\ref{prop: long words in a tuple}, and $5\eta_n =5k^2(2r-1)^{-d(n/2)}$. This concludes the proof since $k<k^2$ and $d(n/2)<n/4$.
\eop

If the size $k$ of the tuple of words is itself a function of $n$, Proposition~\ref{prop: failing ctp} directly yields the following statement.

\begin{corollary}
Let $r=|A|\geq 2$ and let $k(n)$ be an integer function such that $k(n)\leq (2r-1)^{n/2}$.
\begin{enumerate}[(i)]
\item If $k(n)$ is a constant function, then a random $k(n)$-tuple of words in $F(A)$ of length at most $n$ fails the $\log n$-ctp with probability $\O(n^{-\log(2r-1)})$.
\item If $\theta >0$, $k(n)=n^\theta$ and $0<\gamma <1$, then a random $k(n)$-tuple of words in $F(A)$ of length at most $n$ fails the $(2n)^\gamma$-ctp with probability $\O(n^{2\theta} (2r-1)^{-n^\gamma})$.
\item At density $\beta$ (that is, if $k(n)=(2r-1)^{\beta n}$) for $0<\beta <1/8$, and if $4\beta<\gamma<1/2$, then a random $k(n)$-tuple of words in $F(A)$ of length at most $n$ fails the $\gamma n$-ctp with probability $\O((2r-1)^{(4\beta-\gamma)n/2})$.
\end{enumerate}
\end{corollary}

%%%%%%%%%%%%%%%%%%
\section{The uniform membership problem}\label{sec: gmp}

The \emph{Uniform Membership Problem} (UMP) on alphabet $A$ and for an integer $k\geq 1$, is the following: given $w_0$, a word in $F(A)$, and $\vec w=(w_1,\ldots,w_k)$, a $k$-tuple of words in $F(A)$, decide whether $w_0$ belongs to the subgroup $H$ generated by $\vec w$.

The notion of Stallings graphs provides a nice algorithmic solution for the UMP, and to its extension where we also ask for the expression of $w_0$ in a basis of $H$, if $w_0\in H$ \cite{1983:Stallings}:

\paragraph{Algorithm $\algoMP$}
On input a pair $(w_0,\vec w)$ of a reduced word and a $k$-tuple of words in $F(A)$:
\begin{enumerate}[(1)]
\item Compute the Stallings graph $\Gamma(H)$ of $H=\langle \vec w\rangle$.
\item Compute a spanning tree $T$ of $\Gamma(H)$ (which specifies a basis $B$ of $H$, see Section~\ref{sec: subgroups of free groups}).
\item Try reading $w_0$ as a label of a path in $\Gamma(H)$ starting at the root vertex, keeping track of the sequence of edges traversed in the complement of $T$. If one can indeed read $w_0$ in this fashion and the resulting path is a circuit, then $w_0 \in H$ and the sequence of edges not in $T$ yields the reduced expression of $w_0$ in basis $B$ (see Section~\ref{sec: subgroups of free groups}); otherwise $w_0 \not\in H$.
\end{enumerate}

\begin{remark}\label{rk: k is a function of n}
We do not assume the length $k$ of the tuple $\vec w$ to be a constant. We will see it instead as a function of $n=\max|\vec w|$.
\end{remark}

\begin{proposition}\label{prop: worst comp. mp}
The (worst case) complexity of Algorithm $\algoMP$ is $\O(kn\log^*(kn)+rkn+m)$, where $n=\max|\vec w|$, $m=|w_0|$, and $r=|A|$.
\end{proposition}

\proof This is a direct application of Section~\ref{sec: subgroups of free groups}.
\eop

\begin{remark}\label{rk: D without expression}
If we only want to know whether $w_0\in H$, Step~(2) can be skipped, and the complexity is $\O(kn\log^*(kn)+m)$.
\end{remark}

Let us now consider Algorithm $\algoMP$ in more detail, in the case when $\vec w$ has the $d$-ctp for some $d$. In this situation, we can exploit the shape of $\Gamma(H)$ described in Section~\ref{sec: ctp}. In particular, we get a spanning tree by removing one edge from the $\mf_d(w_i)$-labeled path from $\pr_{-i}$ to $\pr_i$ for each $1\le i\le k$, and the corresponding basis is $\vec w$. It does not, actually, matter which edge is removed. As we will see, we do not even need to compute explicitly the full picture of $\Gamma(H)$.

Let $X$ be the $k$-letter alphabet $X=\{x_1,\dots,x_k\}$ and let $\phi\colon F(X)\to F(A)$ be the morphism which maps letter $x_i$ to $w_i$, as in the proof of Proposition~\ref{prop: ctp basis}.

If $w_0 \in H$, then Step ($\algoMP$-3) starts with an \emph{initialization step}, identifying the first letter $x_{i_1}$ of the expression $x_0$ of $w_0$ in basis $\vec w$ and reading $|w_{i_1}|-d$ letters of $w_0$; followed by a potentially \emph{iterated step}, which identifies the next letter in $x_0$. The important observation is that, along each of these steps, a long factor of $w_0$ must be read (of length at least $\min|\vec w|-2d$) and this factor must match one of a fixed collection of at most $2k$ words.

This leads to the family of Algorithms $\algoMPd$ below, indexed by functions $d\colon \mathbb{N}\to \mathbb{N}$, $n\mapsto d(n)$, each of which solves the UMP. We then prove that, for well chosen $d$, the average case complexity of $\algoMPd$ is much lower than the worst case complexity of Algorithm $\algoMP$.

Let $d(n)$ be a non-decreasing function of $n$ such that $d(n) < n/2$.

\paragraph{Algorithm $\algoMPd$}
The input is the $(k+1)$-tuple $(w_0,\dots, w_k)$ of words in $F(A)$. We let $\vec w = (w_1,\dots, w_k)$ and $H=\langle\vec w\rangle$. No assumption (in particular: no ctp assumption) is made about the input. For convenience, we assume that we are also given the $k$-tuple of lengths $(|w_1|,\dots,|w_k|)$ (see Remark~\ref{rk: length of words in input}) and we let $n=\max|\vec w|$. We again let $X=\{x_1,\dots, x_k\}$ and $\phi\colon F(X)\to F(A)$ be given by $\phi(x_i)=w_i$ ($1\leq i\leq k$).
\begin{enumerate}[(1)]
\item Decide whether $\vec w$ has the $d(n)$-ctp and $\min |\vec w|>n/2$. This decision requires computing the set $L(\vec w)$ of length $d(n)$ prefixes of the $w_i$ and $w_i\inv$. This set is recorded in the form of the tree $\Gamma_{d(n)}(\vec w)$, which has at most $2kd(n)$ vertices and edges. There are two cases.
\begin{enumerate}[(a)]
\item If $\vec w$ has the $d(n)$-ctp and $\min |\vec w|>n/2$, go to Step~(2).
\item Otherwise, run Algorithm $\algoMP$ to decide whether $w_0\in H$, and find an expression of $w_0$ in a basis of $H$ if it does.
\end{enumerate}
\item Start reading $w_0$ in $\Gamma_{d(n)}(\vec w)$ from the root vertex. There are two cases.
\begin{enumerate}[(a)]
\item We reach a leaf of $\Gamma_{d(n)}(\vec w)$, say $\pr_{-i}$ ($-k\leq i\leq k$, $i\neq 0$) (necessarily after reading exactly $d(n)$ letters from $w_0$) and the middle factor $\mf_{d(n)}(w_i)$ is a proper prefix of the suffix of $w_0$ starting in position $d(n)+1$ --- that is: $\pr_{-i}\mf_{d(n)}(w_i)$ is a proper prefix of $w_0$. In this case, move the reading head to position $d(n)+|\mf_{d(n)}(w_i)|+1= |w_i|-d(n)+1$ in $w_0$, record $\pr_i$ as the \textsf{last-leaf-visited}, output letter $x_i \in \tilde X$ and go to Step~(3).
\item Otherwise, stop the algorithm and conclude that $w_0 \not\in H$.
\end{enumerate}
\item Suppose that the reading head on $w_0$ is in position $j$, and that the \textsf{last-leaf-visited} is $\pr_i$. Resume reading $w_0$ (from position $j$) in $\Gamma_{d(n)}(\vec w)$, starting at $\pr_i$. There are three cases.
\begin{enumerate}[(a)]
\item We reach the end of $w_0$ while reading it inside $\Gamma_{d(n)}(\vec w)$, landing at the root vertex. In this case, stop the algorithm and conclude that $w_0\in H$.
\item After reading $d'$ letters from $w_0$, we reach a leaf of $\Gamma_{d(n)}(\vec w)$, say $\pr_{-i'}$ (necessarily $2\leq d'\leq 2d(n)$ and $-i' \neq i$), and the word $\mf_{d(n)}(w_{i'})$ is a proper prefix of the suffix of $w_0$ starting in position $j+d'+1$. In this case, move the reading head to position $j+d'+(|w_{i'}| -2d(n))+1$ in $w_0$, record $\pr_{i'}$ as the \textsf{last-leaf-visited}, output letter $x_{i'} \in \tilde X$ and repeat Step~(3).
\item Otherwise, stop the algorithm and conclude that $w_0 \not\in H$.
\end{enumerate}
\end{enumerate}

\begin{theorem}\label{thm: GWP average case}
Let $d(n)$ be a non-decreasing function of $n$ such that $d(n) < n/2$. Algorithm $\algoMPd$ solves the Uniform Membership Problem in $F(A)$ and, if $w_0\in H$, finds an expression of $w_0$ in a basis of $H$.

Let $r=|A|\geq 2$, $0<\delta'<1/4$ and $0<\beta'<\frac{1}{2}-2\delta'$, and suppose that $d(n)\leq \delta'n$. If we restrict the input space to pairs of the form $(w_0,\vec w)$ where $\max|\vec w|=n$ and $\vec w$ is a tuple of length $k\leq (2r-1)^{\beta' n}$, then the average case complexity of $\algoMPd$ is
 $$
\O\left(kd(n)+k^3n(2r-1)^{-d(n/2)}(r+\log^*(kn))+k^2(2r-1)^{-d(n/2)}m\right),
 $$
where $m=|w_0|$. If the space of inputs is further restricted to those inputs where $\vec w$ has the $d(n)$-ctp and $\min|\vec w|>n/2$, then the expected running time is $\O(kd(n))$ --- independent of $|w_0|$.
\end{theorem}

\proof Algorithm $\algoMPd$ always stops because every one of its steps takes a finite amount of time and the only repeated step (Step~(3)) reads a positive number of letters of $w_0$. If $\algoMPd$ stops at Step~(1), then it answers the question whether $w_0 \in H$ and, in the affirmative case, finds an expression for it in a basis of $H$. If it stops at Step~(2), then $w_0 \not\in H$. And if it stops at Step~(3), then $\algoMPd$ outputs a word $x_0$ on alphabet $\tilde X$, one letter at a time (one at the completion of Step~(2) and one at the completion of each iteration of Step~(3) except for the last one). As observed in the description of Step~(3), each new letter cannot be the inverse of the preceding one (because $\Gamma_{d(n)}(\vec w)$ is a tree), so that word $x_0$ is always reduced, that is, $x_0 \in F(X)$. Moreover, the last iteration of Step~(3) concludes either that $w_0\not\in H$, or that $w_0 \in H$ and $x_0=\phi\inv(w_0)$ (that is: $x_0$ is the expression of $w_0$ in basis $\vec w$).

Let us now proceed with bounding the expected running time of $\algoMPd$. We use the following notation:
\begin{eqnarray*}
\mu &=& \min|\vec w|, \\ \probap &=& 2k\,(2r-1)^{-n/2+d(n)}, \\ \probaq &=& 2k\,(2r-1)^{-2-n/2+2d(n)} \enspace=\enspace (2r-1)^{d(n)-2} \probap.
\end{eqnarray*}
Our hypotheses on $k$ and $d$ imply that both $\probap$ and $\probaq$ tend to 0 as $n$ tends to infinity.

Step~(1) first requires comparing the lengths of $w_1,\ldots, w_k$ with $n/2$, deciding whether $\vec w$ has the $d(n)$-ctp and, if so, computing $\Gamma_{d(n)}(\vec w)$. This takes time $\O(kd(n))$ (see Remark~\ref{rk: O(kd)}). If $\vec w$ does not have the $d(n)$-ctp or if $\mu \le n/2$, Step~1 runs Algorithm $\algoMP$, in time $\O(kn\log^*(kn)+ rkn+m)$ (see Proposition~\ref{prop: worst comp. mp}). By Propositions~\ref{prop: long words in a tuple} and~\ref{prop: failing ctp}, this happens with probability $\O(k^2(2r-1)^{-d(n/2)})$ (since we assumed that $d(n)<n/2$ and hence, $k(2r-1)^{-n/2} < k^2(2r-1)^{-d(n/2)}$).

With the complementary probability, $\vec w$ has the $d(n)$-ctp, $\mu>n/2$ and Algorithm $\algoMPd$ proceeds to Step~(2).

We now need to decide in which of the two cases of Step~(2) we are, that is, we need to solve the PPP (Proper Prefix Problem) $2k$ times: for every pair of input words $(u,w_0)$ where $u=\pr_{-i}\mf_{d(n)}(w_i)$ for some $-k\leq i\leq k$, $i\neq 0$. The expected time for this is $\O(k)$ (see Lemma \ref{lm: constant time prefix}). Note that, by the ctp, the output will be positive for at most one of these $u$, thus uniquely identifying the leaf $\pr_{-i}$ of $\Gamma_{d(n)}(\vec w)$ which is first visited when reading $w_0$. Moreover, the probability that the algorithm does not stop here --- and therefore moves to Step~(3) ---, that is, the probability that one of these words $u$ is indeed a proper prefix of $w_0$ is\footnote{The probability that $u$ (of length $\ell$) is a prefix of $w_0$ (of length $m$) is $\O((2r-1)^{-\ell})$. In the present situation, $\ell =|w_i|-d(n)\geq \mu-d(n)$.} $\O(2k(2r-1)^{-\mu+d(n)})$. As $\mu>n/2$, this probability is $\O(\probap)$.

Thus, with probability $\O(\probap)$, we enter a loop where Step~(3) is repeated. Consider one such iteration of Step~(3), starting with the reading head in position $j$ on $w_0$ and vertex $\pr_i$ as the \textsf{last-leaf-visited}. Let $w'_0$ be the suffix of $w_0$ starting at position $j$. To decide in which of the cases of Step~(3) we are, we first consider whether $|w_0|=j+d(n)$, and if so we solve the EqP (Equality problem) on input $(\pr_i\inv,w'_0)$. If indeed $w'_0=\pr_i\inv$, the algorithm stops and concludes that $w_0\in H$. This is done in constant expected time. If $w'_0\neq \pr_i\inv$, we solve the PPP $2k-1$ times, for every input pair $(u,w'_0)$, where $u=\pr_i\inv \pr_{-i'}\, \mf_{d(n)}(w_{i'})$ and $-k \leq i'\leq k$, $i'\neq 0, -i$. By Lemma~\ref{lm: constant time prefix} this is done in expected time $\O(kd(n))$ (the factor $d(n)$ corresponds to the work needed to reduce $\pr_i\inv \pr_{-i'}$ before solving the PPP). The probability that the algorithm continues to a new iteration of Step~(3), namely the probability that one of these words $u$ is a proper prefix of $w'_0$ is $\O((2k-1)(2r-1)^{-(2+\mu-2d(n))})$. Since $\mu>n/2$, we have $(2k-1)(2r-1)^{-(2+\mu-2d(n))} \leq \probaq$.

The expected time required for running Algorithm $\algoMPd$ can be analyzed as follows. Step~(1) runs in expected time
 \begin{align*}
\O(kd(n)) &+ \O(k^2(2r-1)^{-d(n/2)}(kn\log^*(kn)+rkn+m)) \\ &= \O(kd(n) + k^3n(2r-1)^{-d(n/2)}(r+\log^*(kn)) +  k^2(2r-1)^{-d(n/2)}m).
 \end{align*}
With probability $1- \O(k^2(2r-1)^{-d(n/2)})$, Algorithm $\algoMPd$ proceeds to Step~(2).

Step~(2) runs in expected time $\O(k)$. With probability $\O(\probap)$, the algorithm proceeds to Step~(3), and stops with the complementary probability.

Each iteration of Step~(3) runs in expected time $\O(kd(n))$. Step~(3) is repeated with probability $\O(\probaq)$ and the algorithm stops with the complementary probability.

It follows that the expected running time of Step~(2) and the ensuing iterations of Step~(3) is $\O\left( k\left( 1+\probap d(n)(1+\probaq +\probaq^2 +\probaq^3 +\cdots) \right)\right)$, which is
 $$
\O\left(kd(n) \left(1+\frac\probap{1-\probaq}\right)\right).
 $$
Since $\probap$ and $\probaq$ tend to 0, this is $\O(kd(n))$, independently of how many times Step~(3) is iterated. The expected running time of Algorithm $\algoMPd$ is therefore at most
 $$
\O\left(kd(n)+k^3n(2r-1)^{-d(n/2)}\log^*(kn)+k^2(2r-1)^{-d(n/2)}m\right).
 $$
Finally, suppose that the input $(w_0,\vec w)$ is such that $\min|\vec w|> n/2$ and $\vec w$ has the $d(n)$-ctp. Then Step~(1) consists only in computing $\Gamma_{d(n)}(\vec w)$. The expected running time of $\algoMPd$ is therefore, on this smaller set of inputs, $\O(kd(n))$.
\eop

As a corollary, we get upper bounds on the average case complexity of the uniform membership problem.

\begin{corollary}
The Uniform Membership Problem (UMP) for $F(A)$, with input a $k(n)$-tuple of words of length at most $n$, and an additional word of length $m$, can be solved in expected time $C(n,m)$ as follows (where $r=|A|$ is taken to be constant):
\begin{enumerate}[(1)]
\item if $k$ is constant then $C(n,m)=\O(\log n+mn^{-\log(2r-1)})$, improving on its worst case complexity, namely $\O(n\log^*n+m)$.
\item Let $\beta>0$, $0<\gamma<1$. If $k=n^\beta$ then $C(n,m)=\O(n^{\beta+\gamma} +mn^{2\beta}(2r-1)^{-n^{\gamma}})$, improving on its worst case complexity, namely $\O(n^{\beta+1}\log^*n +m)$.
\item If $k=n^\beta$ for some $\beta > 0$, we also have $C(n,m)=\O(n^\beta\log n+ n\log^*n +mn^{-\beta})$.
\item For any $0<\beta<\frac{1}{18}$, if $k=(2r-1)^{\beta n}$ then, for every $0<\epsilon <\frac{1}{8}-\frac{9}{4}\beta$, $C(n,m)=\O(n(2r-1)^{\beta n}+m(2r-1)^{(\frac{9}{4} \beta-\frac{1}{8}+\epsilon)n})$, improving on its worst case complexity, namely $\O(n(2r-1)^{\beta n}\log^*n+m)$.
\end{enumerate}
\end{corollary}

\proof The worst case complexities mentioned in each item follow from Proposition~\ref{prop: worst comp. mp}. For every item, we apply Theorem~\ref{thm: GWP average case} for an appropriate choice of the function $d(n)$:

(1) Suppose that $k$ is a constant function and let $d(n)=\log n$. Note that $(2r-1)^{\log n} = n^{\log(2r-1)}$ and $\lim_{n\to \infty} n^{1-\log(2r-1)}(r+\log^*(kn)) = 0$. Theorem~\ref{thm: GWP average case} then shows that the average case complexity of Algorithm $\algoMPd$ is $\O(\log n +mn^{-\log(2r-1)})$, as announced.

(2) Suppose now that $k=n^\beta$ and let $d(n)=(2n)^\gamma$. We can, again, apply Theorem~\ref{thm: GWP average case}. Since $\lim_{n\to \infty} n^{3\beta+1}(2r-1)^{-n^{\gamma}}(r+\log^*(n^{\beta+1})) = 0$, the average case complexity of $\algoMPd$ is $\O(n^{\beta+\gamma} +mn^{2\beta}(2r-1)^{-n^{\gamma}})$, as announced.

(3) Suppose, again, that $k=n^\beta$ and let $d(n)=\frac{3\beta}{\log(2r-1)}\log(2n)$. Then $(2r-1)^{-d(n/2)}=n^{-3\beta}$. Then Theorem~\ref{thm: GWP average case} shows that, in this case, the average case complexity of $\algoMPd$ is $\O(n^\beta\log n+n\log^*n + mn^{-\beta})$.

(4) Suppose that $k=(2r-1)^{\beta n}$ with $0<\beta <\frac{1}{18}$. This inequality guarantees that $4\beta<\frac{1}{4}-\frac{\beta}{2}$. Let $0<\epsilon <\frac{1}{8}- \frac{9}{4}\beta$ and $\delta =\frac{1}{4}-\frac{\beta}{2}-2\epsilon$ (so that $\delta> 4\beta$), and let $d(n)=\delta n$. Then the hypotheses of Theorem~\ref{thm: GWP average case} are satisfied with $\beta' =\beta$ and $\delta'=\delta$. As a result, the average case complexity of $\algoMPd$ is
 $$
\O\left(n(2r-1)^{\beta n}+n(2r-1)^{(3\beta-\frac{\delta}2)n}\log^*n + m(2r-1)^{(2\beta-\frac{\delta}2)n}\right).
 $$
Since $4\beta<\delta$, we have $3\beta-\frac{\delta}{2}<\beta$ and the second summand is less than the first. Moreover $2\beta-\frac{\delta}{2}=\frac{9}{4}\beta -\frac{1}{8}+ \epsilon$ and the announced result follows.
\eop

%%%%%%%%%%%%%%%%%%
\section{The primitivity and the relative primitivity problems}\label{sec: primitivity}

An element $w$ of a free group $F(A)$ is said to be \emph{primitive} (in $F(A)$) if $F(A)$ admits a basis containing $w$. Equivalently, $w$ is primitive if the cyclic subgroup $\langle w\rangle$ is a free factor of $F(A)$. The \emph{Primitivity Problem} (PrimP) on alphabet $A$ consists in deciding, given a word $w \in F(A)$, whether $w$ is primitive in $F(A)$.

The Primitivity Problem is closely related to the following Whitehead problem: given two words $v, w\in F(A)$, decide whether there exists an automorphism $\phi$ of $F(A)$ such that $\phi(w) = v$. The first step in Whitehead's classical solution to this problem \cite{1936:Whitehead} identifies the minimal length of the automorphic images of $v$ and $w$. The classical solution of PrimP is a by-product of this first step: a word is primitive if and only if its orbit under the action of $\Aut(F(A))$ contains a (and so all) word of length 1. This solution of PrimP is linear in $m=|w|$, but exponential in $r=|A|$ (relying, as it does, on an exploration of the action of the Whitehead automorphisms, whose number is exponential in $r$).

Roig, Ventura, Weil \cite[Fact 3.6]{2007:RoigVenturaWeil} modified Whitehead's algorithm, resulting in an algorithm $\algoPrim$ which solves the primitivity problem in time $\O(m^2r^3)$, where $m=|w|$. We do not describe Algorithm $\algoPrim$ in this paper as we will use it as a black box.

The \emph{Relative Primitivity Problem} (RPrimP) on alphabet $A$ and for an integer $k\geq 1$ is the following: given a word $w_0$ in $F(A)$ and a $k$-tuple $\vec w=(w_1,\dots,w_k)$ of words in $F(A)$, decide whether $w_0$ belongs to $H=\langle \vec w\rangle$ and, if it does, whether it is primitive in $H$.

Solving RPrimP is done naturally by the combination of an algorithm solving the Uniform Membership Problem and, in the case of affirmative answer, computing the expression $x_0$ of $w_0$ in a basis $B$ of $H$, and applying an algorithm for solving PrimP in $F(B)$ (for example, Algorithm $\algoPrim$ mentioned above, with worst-case complexity $\O(m^2r^3)$). By Proposition~\ref{prop: worst comp. mp}, this results in a worst case complexity of $\O(kn\log^*(kn)+rkn+m^2k^3)$ (since the rank of $H$ is at most $k$).

Recently, Shpilrain \cite{2023:Shpilrain} gave an algorithm solving PrimP in $F(A)$ with constant average case complexity. This constant average case complexity assumes, as we have done so far, that the rank $r$ of the ambient free group $F(A)$ is fixed. However, we cannot make this assumption anymore since we need to solve PrimP in free subgroups of $F(A)$, whose rank may be as large as $k$. We therefore revisit this algorithm in detail in Section~\ref{sec: shpilrain} and we recompute its average case complexity to ascertain its dependency in $r$. The average case complexity of the combination of this algorithm with Algorithm $\algoMPd$ (Section~\ref{sec: gmp}) is discussed in Section~\ref{sec: rprim}.

%%%%%%%%%%%%%%%%%%
\subsection{Shpilrain's primitivity algorithm}\label{sec: shpilrain}

Recall that a word $u$ in $F(A)$ is \emph{cyclically reduced} if its last letter is not the inverse of its first letter, that is, if $u^2$ is reduced. It is clear that any word $u$ factors in a unique fashion as $u=vwv\inv$ with $w$ cyclically reduced, and we call $w$ the \emph{cyclic core} of $u$, written $\kappa(u)$. It is immediate that $u$ is primitive if and only if $\kappa(u)$ is.

If $u = x_1\cdots x_n$ is a reduced word of length at least 2, let $W(u)$ be the \emph{Whitehead graph of $u$}, namely the simple (undirected) graph on vertex set $\tilde A$, with an edge from vertex $x$ to vertex $y$ if there exists $1\leq i\leq n$ such that $x_ix_{i+1}=xy\inv$ or $yx\inv$ (here, $x_{n+1}$ stands for $x_1$). Observe that $W(u)$ can be constructed one edge at a time when reading $u$ from left to right, in time $\O(|u|)$.

Recall finally that a vertex $p$ of a connected graph $G$ is a \emph{cut vertex} if deleting $p$ from $G$ (and all the edges adjacent to $p$) results in a disconnected graph. Whitehead showed the following \cite{1936:WhiteheadLMS}.

\begin{proposition}\label{prop: cut vertex}
Let $u$ be a cyclically reduced word of length at least 2 in $F(A)$. If $u$ is primitive, then either $W(u)$ is disconnected, or $W(u)$ admits a cut vertex.
\end{proposition}

Shpilrain's Algorithm $\algoShp$ (slightly modified) is as follows \cite{2023:Shpilrain}. We let
 $$
g(n)=n-\frac1{\log(2r-1)}\, \log(n^4r^6).
 $$

\paragraph{Algorithm $\algoShp$}
On input a reduced word $u\in F(A)$ (given together with its length $n$):
\begin{enumerate}[(1)]
\item Compute $\kappa(u)$, the cyclic core of $u$ --- say, $\kappa(u)=x_1\cdots x_h$. If $h=|\kappa(u)| \leq g(n)$, go to Step~(4). Otherwise, let $i=2$, let $W$ be the graph with vertex set $\tilde A$ and no edges (so that every vertex is its own connected component), and go to Step~(2).
\item Read $x_i$, add the edge $(x_{i-1},x_i\inv)$ to $W$ and update the list of connected components of $W$. If $W$ is connected and has no cut vertex, stop the algorithm: $u$ is not primitive in $F(A)$. Otherwise, if $i<h$, increment $i$ by a unit and repeat Step~(2), and if $i=h$, go to Step~(3).
\item Add to $W$ the edge $(x_h, x_1\inv)$ and update the list of connected components of $W$: if $W$ is connected and has no cut vertex, stop the algorithm: $u$ is not primitive in $F(A)$. Otherwise, go to Step~(4).
\item Run Algorithm $\algoPrim$ on $\kappa(u)$ to decide whether $u$ is primitive in $F(A)$.
\end{enumerate}

Algorithm $\algoShp$ certainly solves the primitivity problem (using Proposition~\ref{prop: cut vertex}), since the graph $W$ constructed in Steps~(2) and~(3) is an increasingly larger fragment of $W(\kappa(u))$. The algorithm stops when either $W$ is connected and has no cut vertex --- in which case $W(\kappa(u))$ has the same property, and $u$ is therefore not primitive in $F(A)$ ---, or when Proposition~\ref{prop: cut vertex} has failed to give us an answer and Algorithm $\algoPrim$ has been called to settle the issue.

Shpilrain \cite{2023:Shpilrain} showed that the average case complexity of Algorithm $\algoShp$ is bounded above by a constant, independent of the length of the input word. This constant does however depend on the ambient rank $r$ and we specify this dependency in Proposition~\ref{prop: Shpilrain's} below. Before we state this proposition, we need to record a few results.

First recall that the number $R_n$ of reduced words of length $n$ in $F(A)$ is $2r(2r-1)^{n-1}$. The number $CR_n$ of cyclically reduced words of length $n$ satisfies
 $$
2r(2r-1)^{n-2}(2r-2) \enspace\leq\enspace CR_n \enspace\leq\enspace 2r(2r-1)^{n-1} \enspace=\enspace R_n.
 $$
It follows that the probability that a reduced word is not cyclically reduced is at most $1-\frac{2r-2}{2r-1}=\frac{1}{2r-1}$ (not exactly $\frac{1}{2r}$, as asserted in~\cite{2023:Shpilrain}, because the first- and last-letter of a reduced word are random variables that are close to but not exactly independent from each other).

The following (obvious!) algorithm computes the cyclic core of a word in $F(A)$.

\paragraph{Algorithm $\algoCR$}
On input a reduced word $u=a_1\cdots a_n$ of length $n$ and, as long as $n\geq 3$: compare $a_n$ with $a_1\inv$; if they are equal, delete the first and last letter of $u$, and repeat this step; if they are different, return the word $u$.

\begin{lemma}\label{lm: complexity algorithm B}
The average case complexity of Algorithm $\algoCR$ is $\O(1)$, independent of the size $r$ of the alphabet.
\end{lemma}

\proof
Let $\probap_n$ be the probability that a length $n$ reduced word is not cyclically reduced. As observed before, $\probap_n \leq \frac{1}{2r-1}$.

Every step of Algorithm $\algoCR$ compares two letters from $\tilde A$, and hence takes time constant time $C$.
On input $u$, of length $n$, Algorithm $\algoCR$ concludes in 1 step (that is the case where $u$ is cyclically reduced) with probability $1-\probap_n$, and otherwise repeats its single step, on a length $n-2$ input. Thus the expected time is bounded above by
 $$
(1+\probap_n +\probap_n\probap_{n-2} +\probap_n\probap_{n-2}\probap_{n-4} +\cdots)\ C \leq \left(\sum_{i\geq 0}(2r-1)^{-i}\right) C.
 $$
Since $\sum_{i\geq 0}(2r-1)^{-i}=\frac{2r-1}{2r-2} \leq \frac{3}{2}$, this concludes the proof.
\eop

We will also use the following fact.

\begin{lemma}\label{lm: length of cyclic core}
The probability that the cyclic core of a length $n$ element of $F(A)$ has length less than or equal to $n-2\ell$ is $\O\left((2r-1)^{-\ell}\right)$.
\end{lemma}

\proof
This is a side product of the proof of Lemma~\ref{lm: complexity algorithm B}: as observed there, Algorithm $\algoCR$ concludes in 1 step with probability $1-\probap_n \leq 1$. It concludes in exactly 2 steps with probability $\probap_{n}(1-\probap_{n-2}) \leq (2r-1)\inv$, and it concludes in $h+1$ steps with probability $\probap_n\probap_{n-2}\cdots\probap_{n-2h+2}(1-\probap_{n-2h}) \leq (2r-1)^{-h}$. Now, $\kappa(u)$ has length $n-2h$ if and only if Algorithm $\algoCR$ terminates in $h+1$ steps. So $|\kappa(u)|\leq n-2\ell$ if $\algoCR$ terminates in at least $\ell+1$ steps, and this happens with probability at most $(2r-1)^{-\ell}\, \sum_i(2r-1)^{-i}$. This quantity is $\frac{2r-1}{2r-2}(2r-1)^{-\ell} \leq \frac32(2r-1)^{-\ell}$.
\eop

An important observation is that, if $u$ is a random word in $F(A)$ of length $n$, then with high probability, $W(u)$ is connected and has no cut vertex. More precisely, the following holds. If $u=x_1\cdots x_n$ is a reduced word of length at least 2, let $W'(u)$ be the simple graph with vertex set $\tilde A$, and with an edge from vertex $x$ to vertex $y$ if there exists $1\leq i<n$ such that $x_i x_{i+1} =xy\inv$ or $yx\inv$. Note that this is almost identical to the definition of the Whitehead graph $W(u)$, except that we do not consider the case where $i=n$. In particular, $W'(u)$ is a subgraph of $W(u)$, and if $W'(u)$ is connected and has no cut vertex, then the same property holds for $W(u)$.

\begin{proposition}\label{prop: proba cut vertex}
Let $r\geq 2$ and let $F=F(A)$, with $|A|=r$. There exists a positive number $\alpha(r)<1-\frac{1}{2}\, r^{-2}$ with the following property: the probability for a word $u$ of length $n$ in $F(A)$ that $W'(u)$ is disconnected, or is connected and has a cut vertex, is $\Theta(\alpha(r)^n)$.
\end{proposition}

\proof
Let $\calG$ be the set of simple graphs on vertex set $\tilde A$ (that is, undirected loop-free graphs without multiple edges). If $G\in \calG$, let $\calA(G)$ be the $\tilde A$-automaton with the same vertex set, whose edges are as follows: for every edge of $G$ connecting vertices $a$ and $b$ ($a,b\in \tilde A$), $\calA(G)$ has a $b\inv$-labeled edge from state $a$ to state $b\inv$ and an $a\inv$-labeled edge from state $b$ to state $a\inv$. Let also $L(G)$ be the set of all words in $\tilde A^*$ which label a path in $\calA(G)$, with no condition on its starting and ending points; in particular, $L(G)$ is a regular language. Finally, let $M(G)$ be the transition matrix of $\calA(G)$, that is, the order $2r$ matrix whose $(a,b)$-entry is 1 if $\calA(G)$ has an edge from vertex $a$ to vertex $b$, and 0 otherwise.

Since $G$ has no loops, every word in $L(G)$ is reduced. In addition, if $u\in F(A)$ is a reduced word, then $u\in L(W'(u))$.

For $G,G'\in \calG$, say that $G\leq G'$ if every edge of $G$ is also an edge of $G'$. Then every labeled edge of $\calA(G)$ is also an edge of $\calA(G')$ with the same label, and we have $L(G) \subseteq L(G')$.

It is clear that if $G\leq G'$ and $G$ is connected and has no cut vertex, then the same holds for $G'$. As a result, if $X$ denotes the set of reduced words $u$ such that $W'(u)$ is disconnected, or is connected but has a cut vertex, then $X\subseteq \bigcup_{i=1}^h L(G_i)$, where $G_1, \ldots, G_h$ are the $\leq$-maximal elements of $\calG$ which are either disconnected, or connected and with a cut vertex.

Let $\lambda_i$ be the growth modulus of $L(G_i)$: then the union of the $L(G_i)$ has growth modulus $\lambda_0 =\max_{1\leq i\leq k}\lambda_i$, and the probability that $u\in \bigcup_{i=1}^h L(G_i)$ is $\Theta(\alpha(r)^n)$, where $\alpha(r)=\frac{\lambda_0}{2r-1}$ (since the growth modulus of the language of all reduced words is $2r-1$).

In order to conclude the proof, we need to establish an explicit upper bound for $\alpha(r)$, as a function of $r$. This is done rather abruptly (following a reasoning similar to that in~\cite{2002:BurilloVentura}): for each $a,b\in \tilde A$ with $b\neq a,a\inv$, let $G_{a,a}$ be obtained from the maximum element of $\calG$ (which has an edge between every pair of distinct vertices in $\tilde A$) by deleting the edge between $a$ and $a\inv$, and $G_{a,b}$ be obtained from the same maximum element by deleting the edge between $a$ and $b$ (see Figure~\ref{fig: transition matrices}).
 \begin{figure}
 $$
\left(\begin{array}{rrrrrrrr}
0 & 0 & 1 & 1 & 1 & 1 & 1 & 1 \\
0 & 0 & 1 & 1 & 1 & 1 & 1 & 1 \\
1 & 1 & 1 & 0 & 1 & 1 & 1 & 1 \\
1 & 1 & 0 & 1 & 1 & 1 & 1 & 1 \\
1 & 1 & 1 & 1 & 1 & 0 & 1 & 1 \\
1 & 1 & 1 & 1 & 0 & 1 & 1 & 1 \\
1 & 1 & 1 & 1 & 1 & 1 & 1 & 0 \\
1 & 1 & 1 & 1 & 1 & 1 & 0 & 1
\end{array}\right)
\textrm{ and }
\left(\begin{array}{rrrrrrrr}
1 & 0 & 1 & 1 & 1 & 1 & 1 & 1 \\
0 & 1 & 0 & 1 & 1 & 1 & 1 & 1 \\
1 & 1 & 1 & 0 & 1 & 1 & 1 & 1 \\
0 & 1 & 0 & 1 & 1 & 1 & 1 & 1 \\
1 & 1 & 1 & 1 & 1 & 0 & 1 & 1 \\
1 & 1 & 1 & 1 & 0 & 1 & 1 & 1 \\
1 & 1 & 1 & 1 & 1 & 1 & 1 & 0 \\
1 & 1 & 1 & 1 & 1 & 1 & 0 & 1
\end{array}\right)
 $$
\caption{The transition matrices of $\calA(G_{a_1,a_1})$ (on the left) and $\calA(G_{a_1,a_2})$ (on the right) for $r = 4$, where the vertex set is $\{a_1,a_{-1},a_2,a_{-2},\dots,a_r,a_{-r}\}$, in that order.}
\label{fig: transition matrices}
 \end{figure}
Then every $G_i$ satisfies either $G_i \le G_{a,a}$ or $G_i\le G_{a,b}$ for some $a,b\in \tilde A$. In particular, $L(G_i)$ is contained in $L(G_{a,a})$ or $L(G_{a,b})$, and hence $\lambda_i$ is less than or equal to the growth modulus of $L(G_{a,a})$ or $L(G_{a,b})$. Since $G_{a,a}$ and $G_{a,b}$ are irreducible and aperiodic, Proposition~\ref{prop: perron-frobenius} shows that these growth moduli are the leading eigenvalues of $M(G_{a,a})$ or $M(G_{a,b})$, respectively.

It should be clear that these growth moduli do not depend on the choice of $a, b \in \tilde A$. Facts~\ref{fact: Gaa} and~\ref{fact: Gab} from the Appendix show that both are at most $(2r-1)(1-\frac{1}{2}\, r^{-2})$, which completes the proof.
\eop

\begin{proposition}\label{prop: Shpilrain's}
The average case complexity of Algorithm $\algoShp$ is $\O\left(\left(\frac{r}{1-\alpha(r)}\right)^2 + r^3\right)$, where $\alpha(r)$ is given by Proposition~\ref{prop: proba cut vertex}, independent of the length of the input word. In particular, this average case complexity is $\O(r^6)$.
\end{proposition}

\proof
Step~(1) of Algorithm $\algoShp$ takes constant expected time, see Lemma~\ref{lm: complexity algorithm B}.

It is a classical result (usually referred to \cite{1973:HopcroftTarjan}) that connectedness and the presence of a cut vertex in a graph with $V$ vertices and $E$ edges can be decided in time $\O(V+E)$. For the graphs occurring in the algorithm, which are subgraphs of the Whitehead graph $W(u)$, we have $V=2r$ and $E\leq 2r(2r-1)$, so $\O(V+E)=\O(r^2)$. It follows that each iteration of Step~(2) takes time $\O(r^2)$, since $W$ has $2r$ vertices, and the same holds for Step~(3).

Finally, Step~(4) takes time $\O(m^2r^3)$, where $m$ is the length of the input word.

By Lemma~\ref{lm: length of cyclic core}, the probability that Step~(1) directly leads to Step~(4), that is, the probability that $|\kappa(u)| \leq g(n)$, is $\O((2r-1)^{-\frac{n-g(n)}2}) =\O(n^{-2}r^{-3})$. So, the contribution of this configuration to the average case complexity of Algorithm $\algoShp$ is $\O(n^{-2}r^{-3}) \O(g(n)^2r^3) = \O(1)$.

Let $q<|\kappa(u)|$ and let $p$ be the length $q-1$ prefix of $\kappa(u)$. Step~(2) is iterated at least $q$ times if the graph $W'(p)$ is disconnected or is connected and has a cut vertex. This happens with probability $\alpha(r)^{q-1}$, where $\alpha(r)$ is given by Proposition~\ref{prop: proba cut vertex}.

Thus the expected running time of Algorithm $\algoShp$ is a big-$\O$ of
 $$
\sum_q \alpha(r)^{q-1}\,q\,r^2 +\alpha(r)^{g(n)}n^2r^3 \enspace\le\enspace \left(\frac{r}{1-\alpha(r)}\right)^2 +\alpha(r)^{g(n)}n^2r^3.
 $$
The inequality above is justified as follows: for $|s|<1$, we have $\sum_q qs^{q-1}=\frac{d}{ds}\left(\sum_q s^q\right)= \frac{d}{ds}\left(\frac1{1-s}\right) =\frac{1}{(1-s)^2}$.

Moreover, we have $\alpha(r)^{g(n)}n^2 < \alpha(r)^nn^2 <\left(1 - \frac1{2r^2}\right)^nn^2$. If $n$ is large enough with respect to $r$, this quantity is less than 1. More precisely, suppose that $\frac {n}{\log n}>4r^2$. Then $\log\left(\alpha(r)^{g(n)}n^2\right) < 2\log n + n\log\left(1 - \frac1{2r^2}\right) < 2\log n - \frac n{2r^2} < 0$. This concludes the proof that the average case complexity of $\algoShp$ is $\O\left(\left(\frac{r}{1-\alpha(r)}\right)^2 + r^3\right)$.

The last part of the statement follows from the observation that $\big( \frac{r}{1-\alpha(r)}\big)^2 <4r^6$.
\eop

\begin{remark}
The Perron-Frobenius theorem can be invoked to show that the spectral radius of $M(G)$ is less than $2r-1$ for each $G\in \calG$ that is not the clique on $\tilde A$. As we saw, we need however an estimate of how much smaller than $2r-1$ these spectral radii are. The method used in the proof of Proposition~\ref{prop: proba cut vertex} is far from optimal: we estimate the spectral radius of $M(G)$ for the graphs obtained from the maximum element of $\calG$ by removing a single edge. Such graphs are far from being disconnected or having a cut vertex. Any upper bound of the spectral radius of the $M(G)$ where $G$ is disconnected or has a cut vertex would lead to an improvement in the expected running time of Algorithm $\algoShp$. There is considerable scope for such an improvement.
\end{remark}

%%%%%%%%%%%%%%%%%%
\subsection{The relative primitivity problem}\label{sec: rprim}

We finally get to the Relative Primitivity Problem, RPrimP: on input a $k$-tuple $\vec w=(w_1,\ldots, w_k)$ of words in $F(A)$ and a word $w_0$ in $F(A)$, along with their lengths, decide whether $w_0$ belongs to $H=\langle \vec w\rangle$ and is primitive in it. As indicated earlier, the idea is essentially to combine an Algorithm $\algoMPd$, for a fast decision of the uniform membership problem, with Shpilrain's Algorithm $\algoShp$, for a fast decision of the primitivity problem, carefully distinguishing  between the situations where $\vec w$ has good properties (the $d(n)$-ctp for a well-chosen function $d$, and the fact that $\min|\vec w|>\frac{1}{2}\max|\vec w|$), which will happen with high probability, and where it does not.

More precisely, consider the following algorithm, parametrized by the choice of a non-decreasing function $d(n)$ such that $d(n)<n/2$.

\paragraph{Algorithm $\algoRPd$}
\begin{enumerate}[(1)]
\item Find out whether $\vec w$ has the $d(n)$-ctp and $\min|\vec w| > n/2$ (this is the first step of Algorithm $\algoMPd$). If one of these properties does not hold, go to Step~(2). If both do, compute $\Gamma_{d(n)}(\vec w)$ and go to Step~(3).
\item Run Algorithm $\algoMP$ on input $(w_0,\vec w)$ to decide whether $w_0\in H=\langle \vec w\rangle$ and, if it does, to compute $x_0$, the expression of $w_0$ in a basis $B$ of $H$. In the latter case, run Algorithm $\algoShp$ on $x_0$ in $F(B)$, to decide whether $x_0$ is primitive in $F(B)$ --- equivalently, whether $w_0$ is primitive in $H$.
\item Run Steps ($\algoMPd$-2) and ($\algoMPd$-3), the latter iterated, to decide whether $w_0\in H$ and, if it does, to compute $x_0$, the expression of $w_0$ in basis $\vec w$. If $w_0\in H$, run Algorithm $\algoShp$ on $x_0$ in the rank $k$ free group $H=\langle \vec w\rangle$.
\end{enumerate}

We prove the following theorem.

\begin{theorem}\label{thm: complexity rprim}
Let $d(n)$ be a non-decreasing function of $n$ such that $d(n)<n/2$. Then Algorithm $\algoRPd$ solves RPrimP.

Let $r=|A|\geq 2$, $0< \delta'<1/4$, $0<\beta'<\frac{1}{2}-2\delta'$ and $\gamma =(2r-1)^{\frac{2\beta'}{1-4\delta'}-1}<1$. Suppose that $d(n)\leq \delta'n$ and $k(n)\leq (2r-1)^{\beta' n}$ for every $n$.
If we restrict the input space to pairs of the form $(w_0, \vec w)$ where $\max|\vec w|=n$ and $\vec w$ is a tuple of length $k(n)$,
then the average case complexity of Algorithm $\algoRPd$ is a big-$\O$ of
 $$
k(n)d(n)+k(n)^2(2r-1)^{-d(n/2)}\left(k(n)n\log^*(k(n)n)+m+k(n)^{6}\right)+ \gamma^{m} k(n)^{6},
 $$
where $m=|w_0|$.

If the input $(w_0,\vec w)$ of RPrimP is limited to those pairs where $\vec w$ has the $d(n)$-ctp and $\min|\vec w|>n/2$, then the average case complexity of $\algoRPd$ is $\O(k(n)d(n)+k(n)^{6} \gamma^{m})$.
\end{theorem}

\proof
It is clear that Algorithm $\algoRPd$ solves RPrimP.

Step~(1) of $\algoRPd$ takes time $\O(k(n)d(n))$, see the analysis of Algorithm $\algoMPd$ in the proof of Theorem~\ref{thm: GWP average case}.

The algorithm moves to Step~(2) with probability $\O(k(n)^2(2r-1)^{-d(n/2)})$, and to Step~(3) with the complementary probability; see Propositions~\ref{prop: long words in a tuple} and~\ref{prop: failing ctp}.

In case we reach Step~(2), as in the proof of Theorem~\ref{thm: GWP average case}, Algorithm $\algoMP$ takes time $\O(k(n)n\log^*(k(n)n)+rk(n)n+m)=\O(k(n)n\log^*(k(n)n)+m)$. If $\algoMP$ concludes that $w_0\in H$, then it also outputs a word $x_0 \in F(B)$, for a certain basis $B$ of $H$; moreover, $|B|\leq k$ and $|x_0| \leq m$. Running Algorithm $\algoShp$ takes time $\O(k(n)^{6})$ in average (see Proposition~\ref{prop: Shpilrain's}).

Otherwise, we reach Step~(3), we are in the situation where $\vec w$ has the $d(n)$-ctp and $\min|\vec w|>n/2$. In particular, the expected running time of Step~$\algoMPd$-2 and all iterations of Step~$\algoMPd$-3 is $\O(k(n)d(n))$, see the proof of Theorem~\ref{thm: GWP average case} (here is where we use the hypothesis on $k(n)$ and $d(n)$).

Moreover, the growth modulus $\gamma_H$ of $H$ is at most $(2k(n)-1)^{\frac 2{n-4d(n)}}$ by Proposition~\ref{prop: growth under ctp}. Then, for every $n$ we have that
 $$
\gamma_H \leq (2k(n)-1)^{\frac {2}{n-4d(n)}}<(2k(n))^{\frac {2}{n-4d(n)}} \leq \left(2\,(2r-1)^{\beta' n}\right)^{\frac {2}{(1-4\delta')n}}=
 $$
 $$
=2^{\frac {2}{(1-4\delta')n}}\cdot (2r-1)^{\frac{2\beta'}{1-4\delta'}}.
 $$
Taking the limit, we get $\gamma_H\leq (2r-1)^{\frac{2\beta'}{1-4\delta'}} <2r-1$ since, by hypothesis, $2\beta'<1-4\delta'$. It follows that the probability that $w_0\in H$ is $\O((\frac{\gamma_H}{2r-1})^m)$ and so, $\O(\gamma^m)$. If indeed $w_0\in H$, we run Algorithm $\algoShp$ on the word $x_0\in F(X)$, in expected time $\O(k(n)^{6})$, see Proposition~\ref{prop: Shpilrain's}.

Thus, the expected running time of this algorithm is bounded above by a big-$\O$ of
 $$
k(n)d(n)+k(n)^2(2r-1)^{-d(n/2)}\left(k(n)n\log^*(k(n)n)+m+k(n)^{6}\right)+ \gamma^mk(n)^{6},
 $$
getting the announced asymptotic estimate.

Finally, if the input $(w_0, \vec w)$ is such that $\vec w$ has the $d(n)$-ctp and $\min|\vec w|>n/2$, then after Step~(1) we go directly to Step~(3) and the expected running time of the algorithm is $\O(k(n)d(n)+\gamma^m k(n)^{6})$.
\eop

As for the uniform membership problem, this gives us upper bounds of the average complexity of the RPrimP for interesting functions $k(n)$.

\begin{corollary}
The relative primitivity problem (RPrimP) for $F(A)$, with input a $k(n)$-tuple of words of length at most $n$, and an additional word of length $m$, can be solved in expected time $C(n,m)$ as follows (where $r = |A|$ is taken to be constant):
 \begin{enumerate}[(1)]
\item If $k$ is constant, then $C(n,m)=\O(\log n+mn^{-\log(2r-1)})$.
\item If $\theta>0$ and $k(n)=n^\theta$, $C(n,m)=\O\left(n^{\theta+\delta} +n^{2\theta}(2r-1)^{-n^\delta}m +n^{6\theta} \left( \frac{2}{2r-1}\right)^m\right)$, for any $0<\delta <1$.
\item If $0<\beta<\frac{1}{58}$ and $k(n)=(2r-1)^{\beta n}$,
 $$
C(n,m)=\O\left(n (2r-1)^{\beta n}+(2r-1)^{-5\beta n}m+(2r-1)^{6\beta n- \frac{1-58\beta}{1-56\beta}m}\right).
 $$
\end{enumerate}
\end{corollary}

\proof
For the case where $k$ is constant, apply Theorem~\ref{thm: complexity rprim} with $d(n)=\log(2n)$ and arbitrary valid values of $\delta'$ and $\beta'$. This shows~(1).

For $k(n)=n^\theta$, choose $d(n)=(2n)^{\delta}$. Then, the hypotheses of Theorem~\ref{thm: complexity rprim} are satisfied with $\delta' =\frac{1}{8}$ and any $\beta'$ such that $0<\beta'<\frac{1}{4}$; the quantity $\gamma$ is then $\gamma =(2r-1)^{4\beta'-1}$. Choosing $\beta'$ such that $(2r-1)^{4\beta'}=2$ yields the announced result. This shows~(2).

Finally, suppose that $k(n)=(2r-1)^{\beta n}$ with $0<\beta <\frac{1}{58}$. Let $\delta =14\beta$ and $d(n)=\delta n$. Note that $\delta <\frac{1}{4}-\frac{\beta}{2}$. Then, the hypotheses of Theorem~\ref{thm: complexity rprim} are satisfied with $\beta' =\beta$ and $\delta' =\delta$. The asymptotic upper bound from that theorem now reads as a big-$\O$ of
 \begin{align*}
(2r-1)^{\beta n}n &+ (2r-1)^{(3\beta-\frac\delta2)n}\log^*n + (2r-1)^{(2\beta-\frac{\delta}{2})n}m \\ & +(2r-1)^{(8\beta-\frac{\delta}{2})n}+ (2r-1)^{6\beta n}\gamma^m,
 \end{align*}
where $\gamma =(2r-1)^{-\frac{1-2\beta-4\delta}{1-4\delta}}$. Since $\delta =14\beta$, we have $\gamma =(2r-1)^{-\frac{1-58\beta}{1-56\beta}}$, the second and fourth summands above are dominated by the first one, and the last summand becomes $(2r-1)^{6\beta n-\frac{1-58\beta}{1-56\beta}m}$. So, our upper bound is a big-$\O$ of
 $$
(2r-1)^{\beta n}n + (2r-1)^{-5\beta n}m + (2r-1)^{6\beta n -\frac{1-58\beta}{1-56\beta}m},
 $$
as announced. This shows~(3), concluding the proof. 
\eop

\section*{Acknowledgments}

The first and second named authors acknowledge partial support from the Spanish Agencia Estatal de Investigaci\'{o}n, through grant PID2021-126851NB-100 (AEI/ FEDER, UE).

%%%%%%%%%%%%%%%%%%%%%%
%%%%%%%%%%%%%%%%%%%%%%
{\small\bibliographystyle{abbrv}
%\bibliography{pwbiblio}}

%
}

%%%%%%%%%%%%%%%%%%%%%%%%
%%%%%%%%%%%%%%%%%%%%%%%%

\newpage
\appendix
\section{Appendix}

\begin{fact}\label{fact: Gaa}
The growth modulus of $L(G_{a,a})$ is
$$\frac12\left(2r-3+\sqrt{(2r+1)^2-8}\right) = (2r-1)\left(1-\frac12 \, r^{-2} - \frac38 \, r^{-4} + \O(r^{-5})\right).$$
\end{fact}

\proof
Let $e_i$ ($-r \le i \le r$, $i\ne 0$) be the column vectors with coordinate at vertex $a_i$ equal to 1 and all other  coordinates equal to 0 --- the standard basis of the dimension $2r$ vector space.

Let $\Maa$ be the transition matrix of $\calA(G_{a,a})$, see Figure~\ref{fig: transition matrices}.

It is an elementary verification that $e_1-e_{-1}$ is in the kernel of $\Maa$, that each $e_i-e_{-i}$ ($2\le i\le r$) is an eigenvector for the eigenvalue 1, and that each $e_i+e_{-i}-e_r-e_{-r}$ ($2\leq i<r$) is an eigenvector for the eigenvalue $-1$. These $2r-2$ vectors, together with $v_1 =e_1+e_{-1}$ and $v_2 =\sum_{i\ge 2}(e_i+e_{-i})$, form a basis of the full space.

Moreover, $\Maa\cdot v_1 = 2 v_2$ and $\Maa\cdot v_2 = (2r-2)v_1 + (2r-3)v_2$. Therefore the other eigenvalues of $\Maa$ are the eigenvalues of $\left(\begin{array}{rr}0 & 2r-2 \\ 2 & 2r-3\end{array}\right)$ and the result follows.
\eop

\begin{fact}\label{fact: Gab}
The growth modulus of $L(G_{a,b})$ is the maximum root of $X^3-(2r-1)X^2+4(r-1)$, and it is bounded above by $(2r-1)(1-\frac12r^{-2})$.
\end{fact}

\proof
Let $e_i$ ($-r \le i \le r$, $i\ne 0$) be as in the proof of Fact~\ref{fact: Gaa} and let $\Mab$ be the transition matrix of $\calA(G_{a,b})$, see Figure~\ref{fig: transition matrices}.

One verifies easily that $e_1-e_2$ is in the kernel of $\Mab$ and that $\Mab \cdot (e_{-1}-e_{-2}) = e_2-e_1$. It also holds that each $e_i - e_{-i}$ ($3\le i\le r$) is an eigenvector for eigenvalue 1, and each $e_i+e_{-i}-e_r-e_{-r}$ ($3\le i < r$) is an eigenvector for eigenvalue $-1$.

These $2r-3$ vectors are linearly independent, and the vectors $v_1 = e_1+e_2$, $v_2 = e_{-1} + e_{-2}$ and $v_3 = \sum_{i=3}^r(e_i + e_{-i})$ complete them to a basis of the full space.

In addition, we have $\Mab \cdot v_1 = 2(v_1+v_3)$, $\Mab \cdot v_2 = v_1 + 2v_2 + 2v_3$ and $\Mab \cdot v_3 = (2r-4)v_1 + (2r-4)v_2 + (2r-5)v_3$.

Thus, in this basis (suitably ordered), the matrix of the linear transformation $\Mab$ consists of the following diagonal blocks:
$$\left(\begin{array}{rr} 0 & -1 \\ 0 & 0\end{array}\right), \quad (1) \textrm{ ($r-2$ times),}\quad (-1) \textrm{ ($r-3$ times),} \quad \left(\begin{array}{rrr} 2 & 1 & 2r-4 \\ 0 & 2 & 2r-4 \\ 2 & 2 & 2r-5\end{array}\right).$$
In particular, the remaining eigenvalues are the roots of the characteristic polynomial of that $(3, 3)$-matrix, namely $P(X) = X^3-(2r-1)X^2+4(r-1)$. We note that the local extrema of $P(X)$ are at $0$ and $\frac23(2r-1)$, $P(0) = 4(r-1) >0$ and $P(r) = -r^3 + r^2 + 4r - 4$, which is negative for all $r \ge 2$. Therefore $P(X)$ has three real roots.

Since $P(2r-1) = 4(r-1) > 0$, the leading eigenvalue of $\Mab$ sits between $\frac23(2r-1)$ and $2r-1$. For a closer estimate, let $\delta = (2r-1)(1-\frac12\,r^{-2})$. Then
\begin{align*}
P(\delta) &= 4(r-1) + (2r-1)^3\left(1-\frac1{2r^2}\right)^2\left(1-\frac1{2r^2} - 1\right) \\
&= 4(r-1) - \frac{(2r^2-1)^2}{8r^6}(2r-1)^3 \\
& = \frac{16 r^{6} + 8 r^{5} - 44 r^{4} + 16 r^{3} + 8 r^{2} - 6 r + 1}{8r^6},
\end{align*}
which is positive when $r\ge 2$. Thus the leading eigenvalue of $\Mab$, which is the growth modulus of
$L(G_{a,b})$, is at most $\delta = (2r-1)(1-\frac12\,r^{-2})$ as announced.
\eop

\begin{remark}
Facts~\ref{fact: Gaa} and~\ref{fact: Gab}, while mathematically elementary, would have been very difficult to establish without the help of a versatile computer algebra system. The authors are grateful to the developers of SageMath \cite{2022:The-Sage-Developers}.
\end{remark}

%%%%%%%%%%%%%%%%%%%%%%%%
\end{document}